\documentclass{amsart}
\usepackage[utf8]{inputenc}

\usepackage{tikz}
\usepackage{xcolor}

\usepackage{amsmath}
\usepackage{amsthm}
\usepackage{amssymb}

\usepackage{booktabs}

\usepackage[a4paper]{geometry}

\usepackage{hyperref}

\newtheorem{thm}{Theorem}[section]

\newtheorem*{ack}{Acknowledgements}

\theoremstyle{definition}

\newcommand{\bbC}{\mathbb{C}}
\newcommand{\bbR}{\mathbb{R}}

\newcommand{\bbZ}{\mathbb{Z}}

\newcommand{\bbP}{\mathbb{P}}

\DeclareMathOperator{\diag}{diag}
\DeclareMathOperator{\SL}{SL}
\DeclareMathOperator{\Bl}{Bl}

\title{CscK metrics on rank one spherical Fano fourfolds}

\author{Thibaut Delcroix}
\address{Univ Montpellier,
CNRS, Montpellier, France}
\email{thibaut.delcroix@umontpellier.fr}

\subjclass[2020]{}
\keywords{}
\date{\today}

\begin{document}

\begin{abstract}
This four-pages note is an invitation to explore explicit K-stability for arbitrary Kähler classes of low dimension and low rank spherical varieties. 
We apply our simple combinatorial criterion of K-stability of rank one spherical varieties to the example of the blowup of the product of two copies of the projective planes along the diagonal, and obtain strong indication that it admits cscK metrics in every Kähler classes. 
\end{abstract}

\maketitle

\section{Introduction}

The purpose of this note is to further illustrate how the author's result in \cite{Delcroix_RK1} allows to check for the existence of cscK metrics on rank one spherical manifolds. 
In view of the classification achieved in \cite{Delcroix-Montagard}, in this note we almost settle the question of existence of cscK metrics in arbitrary Kähler classes on rank one spherical Fano fourfolds. 

Recall that a normal complex algebraic variety \(X\) equipped with a regular action of a connected reductive group \(G\) is spherical if any Borel subgroup \(B\) of \(G\) acts on \(X\) with an open orbit. 
A fundamental invariant of such a \(G\)-variety is its rank, defined as the rank of the \(\bbZ\)-module generated by the \(B\)-weights of \(B\)-semi-invariant rational functions on \(X\). 
In \cite{Delcroix-Montagard}, together with Pierre-Louis Montagard, we classified the locally factorial Fano spherical \(G\)-varieties of dimension four or less, and rank two or less. 
We further checked, using \cite{Delcroix_KSSV}, which of these Fano varieties admit Kähler-Einstein metrics. 

Kähler-Einstein metrics are only one among several possible candidates for canonical Kähler metrics, and they are associated with the choice of the anticanonical polarization. 
For a general polarization, the canonical metrics of choice are usually extremal Kähler metrics as defined by Calabi \cite{Calabi_1982}, which, when the Futaki invariant vanishes, are constant scalar curvature Kähler (cscK) metrics. 
The Yau-Tian-Donaldson conjecture states that existence of extremal Kähler metrics on a polarized manifold \((X,L)\) should be equivalent to (\(G\)-uniform) K-stability. 
This statement was proved in the particular case of smooth spherical varieties and for cscK metrics in the appendix of \cite{Delcroix_KSSV2} by Odaka, relying on the work of Li \cite{Li_2022}. In \cite{Delcroix_RK1} we further proved that the K-stability condition for rank one spherical manifolds consists in a single combinatorial condition. 

The purpose of this note is to write the condition explicitly for the example of the blowup \(X\) of \(\bbP^2\times \bbP^2\) along the diagonal. 
By our classification in \cite{Delcroix-Montagard}, a smooth spherical rank one Fano manifold of dimension four or less is either this \(X\), or horospherical (under a possibly different subgroup). 
It is strongly expected (and would follow from the proof of \cite{Delcroix_RK1}, adapted to extremal Kähler metrics) that on a horospherical rank one manifold, all Kähler classes admit extremal Kähler metrics, so that on such a manifold, a Kähler class contains a cscK metric iff it has vanishing Futaki invariant. In view of these observation, dealing with the case of  \(\Bl_{\diag \bbP^2}(\bbP^2\times \bbP^2)\) thus essentially puts an end to the quest for cscK metrics on spherical rank one Fano fourfolds. 
Once the condition is written out, any numerical analysis program gives evidence that \(\Bl_{\diag \bbP^2}(\bbP^2\times \bbP^2)\) admits cscK metrics in all Kähler classes, but it is difficult to work out a short formal proof of this. 
We instead examine the condition near the walls of the ample cone, as well as on a remarkable subset, and show that it is satisfied in these cases (and close to these cases, since the condition is open). 

\section{The manifold as a spherical variety, and its ample cone}

We consider the group \(\SL_3(\bbC)\). 
We consider the Borel subgroup \(B\) of upper triangular matrices, and the maximal torus \(T\) of diagonal matrices. 
To this data is also associated a root system \(\Phi\) and a positive root system \(\Phi^+\). 
We denote by \(\alpha_1\) and \(\alpha_2\) the simple roots, so that \(\Phi^+=\{\alpha_1,\alpha_2, \alpha_1+\alpha_2\}\). 
Let \(\varpi_1\) and \(\varpi_2\) denote the corresponding fundamental weights. 
For use in the next section, let us record that the half-sum of positive roots is 
\( \varpi = \alpha_1+\alpha_2=  \varpi_1+\varpi_2 \), 
that \(\alpha_2=2\varpi_2-\varpi_1\), and that we have
\[ \frac{\langle \alpha_1, x_1\varpi_1+x_2\varpi_2\rangle}{\langle \alpha_1, \varpi_1+\varpi_2\rangle} = x_1 \qquad \frac{\langle \alpha_2, x_1\varpi_1+x_2\varpi_2\rangle}{\langle \alpha_1, \varpi_1+\varpi_2\rangle} = x_2 \qquad \frac{\langle \alpha_1+\alpha_2, x_1\varpi_1+x_2\varpi_2\rangle}{\langle \alpha_1, \varpi_1+\varpi_2\rangle} = \frac{x_1+x_2}{2} \]

Consider the manifold \(\bbP^2\times \bbP^2\) equipped with the diagonal action of \(\SL_3(\bbC)\). 
It is spherical, as well as its blowup \(X\) along the diagonal embedding of \(\bbP^2\) which is a codimension 2 orbit of \(\SL_3(\bbC)\). 

We recall the combinatorial data associated to \(X\) from \cite[\textbf{4-1-3}, Section~7.2]{Delcroix-Montagard}. 
The weight lattice \(M\) is generated by \(\alpha_2\), and the unique spherical root is \(\alpha_2\). 
We identify an element \(\nu\) of \(N\) with the image of \(\alpha_2\) by \(\nu\) in \(\bbZ\), thus we identify \(N\) with \(\bbZ\).  
There are three colors \(\mathcal{D}=\{\clubsuit,\heartsuit,\diamondsuit\}\), and one \(G\)-stable prime divisor \(E\) (the exceptional divisor of the blowup) with images in \(N\) given by \(\clubsuit|_M=E|_M=-1\) and \(\heartsuit|_M=\diamondsuit|_M=1\).
Note that the divisor \(b\diamondsuit +c\heartsuit\) represents the pull-back of the line bundle \(O(b,c)\) on \(\bbP^2\times \bbP^2\) under the blowup map. 
By \cite{Brion_1989}, the Picard group of \(X\) is generated by these divisors, and the unique relation is given by \(\diamondsuit + \heartsuit = E+\clubsuit\). As a consequence, any Cartier divisor may be written as 
\[ D(a,b,c) = -aE+b\diamondsuit +c\heartsuit \]
and the corresponding polytope is the segment 
\[ \Delta(a,b,c) = [\sup(-b,-c)\alpha_2,\inf(-a,0)\alpha_2] \subset M\otimes \bbR\]

Since \(\SL_3\) has no non-trivial characters, a \(G\)-linearized line bundle \(L\) has a unique linearization. 
There exists a (unique up to scaling) %up to passing to a multiple of \(L\)? (fo existence...) in these two examples, must be unnecessary
\(B\)-equivariant rational section \(s\) of \(L\) whose divisor is of the form \(D(a,b,c)\), and the \(B\)-weight of this section is \((b+c)\varpi_2\). 
Thus the moment polytope of \(L\) is 
\[ \Delta^+(a,b,c) = (b+c)\varpi_2 + [\sup(-b,-c)\alpha_2,\inf(-a,0)\alpha_2] \subset X^*(B)\otimes \bbR\]

By Brion's ampleness criterion \cite{Brion_1989}, the line bundle \(L\) is ample if and only if \(0< a< \inf(b,c)\). 
Since \(b\) and \(c\) play a completely symmetric role, we may as well assume \(b\leq c\), and up to passing to real line bundles, we may assume up to scaling that \(b+c=1\) if \(L\) is ample. 
Thus the set of ample real line bundles up to scaling and symmetry is parametrized by the triangle 
\( \{(a,b)\mid 0<a<b \leq \frac{1}{2}\}\).
For an ample real line bundles \(L(a,b,1-b)\) with \(0<a<b\leq \frac{1}{2}\), its moment polytope is 
\[\Delta^+ = (b+c)\varpi_2 + [-b,-a]\alpha_2 \]

\section{The K-stability condition for \(L(a,b,c)\)}

By \cite{Delcroix_RK1}, there is a single combinatorial condition to check. 
We use the notations from \cite[Section~2]{Delcroix_RK1} to write it explicitly. 
In these notations, we have \(R_X^+=\Phi^+\), \(\sigma=\alpha_2\), \(\chi=\varpi_2\),  \(s_-=-b\), \(s_+=-a\), 
\[ P(t) = \frac{-2t^3-3t^2-t}{2} \]
and 
\[ Q(t) = \frac{1-3t^2}{2} \]
Let 
\begin{align*} 
C=C(a,b) := & 
\left(\int_{-b}^{-a}P(t)\mathop{dt}\right)\left(-bP(-b)-aP(-a)+2\int_{-b}^{-a}tQ(t)\mathop{dt}\right) \\
 & 
 - \left(\int_{-b}^{-a}tP(t)\mathop{dt}\right)\left(P(-b)+P(-a)+2\int_{-b}^{-a}Q(t)\mathop{dt}\right)
\end{align*}
The K-stability condition for the real ample line bundle \(L(a,b,1-b)\) is 
\[ C(a,b) > 0  \]
 
We compute 
\[ \int_{-b}^{-a}P(t)\mathop{dt}=\frac{b^4-a^4-2(b^3-a^3)+b^2-a^2}{4} \]
\[ \int_{-b}^{-a}tP(t)\mathop{dt}=\frac{-24(b^5-a^5)+45(b^4-a^4)-20(b^3-a^3)}{120} \]
\[ P(-b)+P(-a)+2\int_{-b}^{-a}Q(t)\mathop{dt} = \frac{4a^3-3(a^2+b^2)+3b-a}{2} \]
\[ -bP(-b)-aP(-a)+2\int_{-b}^{-a}tQ(t)\mathop{dt} = \frac{-b^4-7a^4+6(b^3+a^3)+4a^2}{4} \]
yielding an unscrutable expression for \(C(a,b)\), although one can very easily check, given explicit values of \(a\) and \(b\), whether \(C(a,b)>0\). 
Numerical computations indicate that \(C(a,b)>0\) for \(0<a<b\leq \frac{1}{2}\), so that \(X\) admits cscK metrics in every Kähler classes, but there does not seem to be a simple argument to prove this. 
Instead, we verify it on the boundary cases. 

Assume first that \(b=\frac{1}{2}\). This segment is actually not on the boundary of the ample cone, since we used symmetry in \(b\) and \(c\). Actually, a multiple of the anticanonical line bundle lies in this segment, and the anticanonical line bundle was proven to be K-stable in \cite{Delcroix-Montagard}. 
We have, after factorization,
\begin{align*}
C(a,\frac{1}{2}) & = \frac{(1-2a)^2}{61440}(576a^6-2496a^5-464a^4+1888a^3+524a^2+884a+249) \\
& > \frac{(1-2a)^2}{61440}(576a^6+1032a^3+524a^2+884a+249) \\
\intertext{by using \(a<\frac{1}{2}\), hence}
C(a,\frac{1}{2}) & > \frac{249(1-2a)^2}{61440} >0
\end{align*}
by using \(a>0\). 
Since the condition is open, any ample line bundle \(L(a,b,c)\) with \(\frac{c-b}{b+c}\) small enough admits a cscK metric. 

Assume now that \(a=b\). 
Obviously, \(C(a,b)=0\), but we can factor out a \((b-a)\) to check the sign as one approaches this case. 
We have 
\[ \left.\frac{C(a,b)}{b-a}\right|_{b=a} = a^3(1-2a)(1-a) > 0 \]
As a consequence, any ample line bundle \(L(a,b,c)\) with \(\frac{b-a}{b+c}\) small enough admits a cscK metric. 

Finally, we consider the case \(a=0\). 
This corresponds to semi-ample classes that are pulled-back from ample line bundles on \(\bbP^2\times \bbP^2\), and the condition reduces to the condition we would get for these bundles by considering \(\bbP^2\times \bbP^2\) as a rank one \(\SL_3(\bbC)\)-spherical variety. 
Since any Kähler class on \(\bbP^2\times \bbP^2\) admits a cscK metric, the condition must be satisfied. 
We check this directly:
\begin{align*} 
C(0,b) & = \frac{b^4(1-b)(5b^3-11b^2-15b+20)}{80} \\
& > \frac{b^4(1-b)\frac{39}{4}}{80} > 0 
\end{align*}
by using \(0<b<\frac{1}{2}\). 
As a consequence, any ample line bundle \(L(a,b,c)\) with \(\frac{a}{b+c}\) small enough admits a cscK metric. 
To summarize:

\begin{thm}
Consider the blowup \(\pi: X=\Bl_{\diag \bbP^2}\bbP^2\times \bbP^2\to \bbP^2\times \bbP^2\), with the exceptional divisor \(E\). 
Let \(0<a<b\leq c\) be integers, and consider the ample line bundle \(L(a,b,c) = \pi^*\mathcal{O}_{\bbP^2\times \bbP^2}(b,c) \otimes \mathcal{O}(-aE)\). 
There exists \(\varepsilon>0\) such that if \(\inf(\frac{a}{b+c},\frac{b-a}{b+c},\frac{c-b}{b+c})<\varepsilon\), then there exists a cscK metric in \(c_1(L(a,b,c))\).  
From numerical investigations, \(\varepsilon = +\infty\) should also work. 
\end{thm}

\begin{ack}
The author is partially funded by ANR-21-CE40-0011 JCJC project MARGE.
\end{ack}

\end{document}